\documentclass{amsart}
\usepackage{graphicx}
\usepackage{amsmath, amssymb, amsfonts}
\usepackage{hyperref}
\usepackage[shortalphabetic]{amsrefs}
\usepackage{mathptmx}    

\begin{document}
\title{Contracting self-similar solutions of nonhomogeneous curvature flows} 

\author[J. A. McCoy]{James A. McCoy}
\address{Priority Research Centre Computer Assisted Research Mathematics and Applications, School of Mathematical and Physical Sciences, University of Newcastle, University Drive, Callaghan, NSW 2308, Australia}
\email{James.McCoy@newcastle.edu.au}
\keywords{curvature flow, parabolic partial differential equation, hypersurface}
\subjclass[2000]{53C44}

\begin{abstract}
A recent article \cite{LL20} considered fully nonlinear contraction of convex hypersurfaces by certain nonhomogeneous functions of curvature, showing convergence to points in finite time in cases where the speed is a function of a degree-one homogeneous, concave and inverse concave function of the principle curvatures.  In this article we consider self-similar solutions to these and related curvature flows that are not homogeneous in the principle curvatures, finding various situations where curvature pinched, convex, mean-convex or even general closed hypersurfaces contracting self-similarly are necessarily spheres.  
\end{abstract}

\keywords{curvature flow, parabolic partial differential equation, self-similar solution}

 \thanks{This research was supported by Discovery Grant DP180100431 of the Australian Research Council.  Part of this work was completed while the author was at the Okinawa Institute of Science and Technology as part of the Visiting Mathematics Professors program.  The author is grateful for these sources of support.}

\maketitle

\section{Introduction} \label{S:intro}
In \cite{M11}, the author considered contracting self-similar solutions of fully nonlinear curvature contraction flows whose speeds were homogeneous functions of the principle curvatures.  This work extended earlier results of Huisken for the mean curvature flow \cite{H90}, where it was shown that a compact hypersurface with nonnegative mean curvature contracting self-similarly under the mean curvature flow is necessarily a sphere.  In the case of surfaces of dimension $2$, the condition of nonnegative mean curvature was replaced by the requirements that the surface be embedded and have genus $0$ by Brendle \cite{B16}.

For a survey of more general self-similar shrinkers of the mean curvature flow we refer to \cite{DLN18}.  For other flows, we refer the reader to the discussion in \cite{M11}, noting that few results are available apart from those for flows by powers of the Gauss curvature.  Let us here briefly update recent developments.  In \cite{MMW15} the author, together with Mofarreh and V-M Wheeler extended the results of \cite{Asurf} for contracting surfaces to the case of axially symmetric hypersurfaces., including  corresponding characterisation results for when axially symmetric hypersurfaces contracting self-similarly are necessarily spheres.  In \cite{BCD17}, the authors proved that closed, strictly convex hypersurfaces contracting self-similarly by powers $\alpha \geq \frac{1}{n+2}$ of the Gauss curvature are spheres if the inequality is strict, or ellipses in the equality case (see also \cite{CD16} for some powers).  This removes the extra conditions of \cites{M11, MMW15} required to conclude the hypersurfaces are spheres in this case.  The result was also shown slightly earlier with an additional symmetry assumption in \cite{AGN16}.  

The result of \cite{BCD17} was generalised to powers $\alpha \geq \frac{1}{k}$ of the elementary symmetric functions $\sigma_k$ of the principle curvatures, $1\leq k \leq n-1$ in \cite{GLM18}, showing closed, strictly convex hypersurfaces contracting self-similarly must be round spheres.  Further generalisations are in \cite{C19} to speeds $\left( \frac{\sigma_n}{\sigma_{n-k}}\right)^{\alpha}$ for $\alpha>\frac{1}{k}$ and in \cite{CG19} to speeds $\left( \frac{\sigma_k}{\sigma_{\ell}} \right)^{\alpha}$ for $0\leq \ell < k < n$ and $\alpha>\frac{1}{k-\ell}$.  There are also some results for closed self-similar hypersurfaces in other ambient spaces (eg \cite{GM19}) and for higher order curvature flows such as the curve diffusion flow \cite{E15}.  We will not discuss these settings further here, nor other interesting related problems including other types of self-similarity and the case of self-similar evolving curves.

The structure of this article is as follows.  In Section \ref{S:prelim} we set up the problem, state our main result and describe the structure conditions on the flow speed.  We give some example speeds that fit our requirements, we provide some geometric estimates that are needed in the later sections and we give the general structure of the geometric partial differential equation that is at the heart of establishing our results.  In the subsequent sections we provide proofs of the main result in each of the cases of particular conditions on the speed.

\section{Preliminaries} \label{S:prelim}
\newtheorem{Fconditions}{Conditions}[section]
\newtheorem{Phiconditions}[Fconditions]{Conditions}
\newtheorem{main}[Fconditions]{Theorem}
\newtheorem{eeqns}[Fconditions]{Lemma}
\newtheorem{Consign}[Fconditions]{Lemma}
\newtheorem{Fconvprops}[Fconditions]{Lemma}
\newtheorem{curvest}[Fconditions]{Lemma}

Let $M_{0}$ be a compact, convex hypersurface of dimension $n\geq 2$, without boundary, smoothly embedded in $\mathbb{R}^{n+1}$ and represented  by some diffeomorphism $X_{0}:\mathbb{S}^{n} \rightarrow X_{0}\left( \mathbb{S}^n \right)= M_{0}\subset \mathbb{R}^{n+1}$.  We consider the family of maps $X_{t}=X \left( \cdot, t \right)$ evolving according to 
\begin{equation} \label{E:theflow}
\begin{split}
  \frac{\partial}{\partial t}X \left( x,t \right) 
= - \Phi\left( F\left( \mathcal{W} \left( x,t\right) \right) \right) &\nu\left( x,t \right) \mbox{ } x\in \mathbb{S}^n,\mbox{ }  0 < t\leq T\leq \infty \\
X \left( \cdot, 0\right) 
&= X_{0} \mbox{,}
\end{split}
\end{equation}
where  $\mathcal{W}\left( x, t\right) $ is the matrix of the Weingarten map of $M_{t}=X_{t}\left( \mathbb{S}^n \right)$ at the point $X_{t}\left( x\right)$ and $\nu\left( x, t\right)$ is the outer unit normal to $M_{t}$ at $X_{t}\left( x\right)$.  

We will shortly describe the properties of the functions $F$ and $\Phi$ in detail but for now let us first note that $\Phi$ is a positive function so \eqref{E:theflow} is a contraction flow.  In this article we are specifically interested in hypersurfaces that contract self-similarly, that is
\begin{equation} \label{E:ssc}
  X\left( x, t\right) = \psi\left( t\right) X\left( x, 0 \right) \mbox{,}
\end{equation}
for some nonnegative, monotone decreasing function $\psi$ with $\psi\left( 0\right) = 1$.  As $\Phi$ and $F$ will have properties ensuring that \eqref{E:theflow}   parabolic, exterior and interior spheres evolving coincidently with $X$ remain disjoint to $M_t$ and provide estimates above and below on the maximal existence time $T$.  In our setting the exterior shrinking sphere gives an estimate from above on the time by which the function $\psi\left( t\right)$ has decreased to zero.

To our knowledge, previous work on self-similar contracting hypersurfaces has always used `separation of variables' to solve an ordinary differential equation for $\psi$ and then work with the corresponding elliptic equations for geometric quantities associated $M_0$ to deduce the potential shapes that contract self-similarly.  This relies crucially on homogeneity of the speed.  As our speeds generally are not homogeneous we cannot perform separation of variables so we have to keep all factors of $\psi$ and its derivative in our equations.  Specifically, differentiating \eqref{E:ssc} with respect to $t$
$$\frac{\partial X}{\partial t} \left( x, t\right) = \psi'\left( t\right) X\left( x, 0 \right) \mbox{,}$$
where we will use $'$ to denote the derivative of a function of one variable.  Substituting this into \eqref{E:theflow} we obtain
$$\psi'\left( t\right) X\left( x, 0 \right) =  - \Phi\left( F\left(  \frac{1}{\psi\left( t\right)}\mathcal{W} \left( x,0\right) \right) \right) \nu\left( x,0 \right)$$
where we have used that, in view of \eqref{E:ssc}
$$\mathcal{W}\left( x, t\right) = \frac{1}{\psi\left( t\right)} \mathcal{W}\left( x, 0\right) \mbox{.}$$
and $\nu\left( x, t\right) = \nu\left( x, 0\right)$.  Taking now the inner product with $\nu\left( x, 0 \right)$ we obtain our equation for self-similar hypersurfaces in this setting:
\begin{equation} \label{E:sse}
  \left< X\left( x, 0 \right), \nu\left( x, 0 \right) \right> = -\frac{1}{\psi'\left( t\right)} \Phi\left( F\left(  \frac{1}{\psi\left( t\right)}\mathcal{W} \left( x,0\right) \right) \right) \mbox{,}
  \end{equation}
  where we have divided through by $\psi'\left( t\right)$ which is a positive function until the extinction time $T$.

Now we describe the properties of the functions $F$ and $\Phi$.  We denote by $\Gamma_+$ the positive cone, $\Gamma_+=\left\{ \kappa = \left( \kappa_1, \ldots, \kappa_n \right) : \kappa_i >0 \mbox{ for all }i=1, \ldots, n \right\}$.  

The function $F$ should have the following properties:

\begin{Fconditions} \label{T:Fconds}
\mbox{}
\begin{enumerate}
  \item[\textnormal{a)}] $F\left( \mathcal{W} \right) = f\left( \kappa \left( \mathcal{W} \right) \right)$ where $\kappa\left( \mathcal{W}\right)$ gives the eigenvalues of $\mathcal{W}$ and $f$ is a smooth, symmetric function defined on an open, symmetric cone $\Gamma \supseteq \Gamma_+$. 
  \item[\textnormal{b)}] $f$ is strictly increasing in each argument: $\frac{\partial f}{\partial \kappa_{i}} > 0$ for each $i=1, \ldots, n$ at every point of $\Gamma$.
  \item[\textnormal{c)}] $f$ is degree-one homogeneous: $f\left( k \kappa \right) = k f\left( \kappa\right)$ for any $k>0$.
  \item[\textnormal{d)}] $f$ is strictly positive on $\Gamma$ and $f\left( 1, \ldots, 1\right) = 1$.
\end{enumerate}
\end{Fconditions}

\noindent \textbf{Remarks:}
\begin{enumerate}
  \item Examples of $F$ and their corresponding cones $\Gamma$ are given in \cite{M11}.  In this article we do require $M$ to be convex (to be able to establish signs on the additional terms that arise in equations for geometric quantities); the cones on which $f$ are defined are generally pinching cones of the form
  $$\Gamma_\varepsilon = \left\{ \kappa \in \Gamma^+: \kappa_i \geq \varepsilon \kappa_j \right\} \mbox{.}$$
  When $\varepsilon = 1$ all principle curvatures are equal and only spheres are possible.  Smaller $\varepsilon$ corresponds to weaker curvature pinching.    
  \item Often Conditions \ref{T:Fconds} are accompanied by an additional second derivative condition on $f$.  We will prove our result under several sets of conditions, and mention a second derivative condition on $f$ where it is needed.
  \item In the case $\Phi\left( F \right) = F^\alpha$, $\alpha \geq 1$, self-similar solutions of \eqref{E:theflow} were considered in \cite{M11} (with the axially symmetric hypersurface case considered in \cite{MMW15}).  More generally contraction flows of convex hypersurfaces with speeds $F$ and $F^\alpha$ for $\alpha>0$ have been widely considered; we refer the reader to \cite{AMZ11} for example and the references contained therein for general behaviour of such flows.  In a forthcoming article \cite{M20} we will consider the case of general nonhomogeneous curvature contraction, complementary to \cite{LL20}.
  \end{enumerate}

For our results we will require the function $\Phi: \left[ 0,\infty \right)\rightarrow \mathbb{R}$ to be at least twice differentiable.  For our various results we will require Conditions \ref{T:Phiconditions} a) to c) below; some will require also d).  Our requirements do not rule out the possibility that results might be possible with different requirements on $\Phi$ (or indeed requirements on $F$) using different approaches.

\begin{Phiconditions} \label{T:Phiconditions}
\begin{itemize}
  \item[\textnormal{a)}] $\Phi\left( 0 \right) = 0$;
  \item[\textnormal{b)}] $\Phi'\left( z\right) > 0$ for all $z>0$; 
  \item[\textnormal{c)}] $\Phi''\left( z\right) \geq 0$ for all $z>0$, 
    \item[\textnormal{d)}] There is a constant $c>0$ such that, for all $z>0$, $z \left| \Phi''\left( z\right) \right| \leq c\, \Phi'\left( z\right)$ for some constant $c>0$; 
\end{itemize}
\end{Phiconditions}

In particular, Condition b) is needed for parabolicity of \eqref{E:theflow}; this is an essential requirement.   Condition c), which implies also that $\Phi'\left( z\right) z - \Phi \left( z\right) \geq 0$ for all $z>0$, is needed to ensure certain terms in our equations have the correct sign for applying the maximum principle.  Condition d) is used in some cases to allow an unsigned derivative term to be absorbed into a good term provided $M$ is sufficiently curvature-pinched.\\

\noindent \textbf{Remarks:}
\begin{enumerate}
  \item The above conditions may be compared with some of those used in \cites{BP, BS, LL20} where constrained curvature contraction flows were considered.  In most cases our conditions here are a little less restrictive than in \cite{M20}.\\
  \item Some example functions $\Phi$ include
    \begin{enumerate}
  \item $\Phi\left( z\right) = \sum_{i=1}^{\ell} c_i z^{k_i}$, for constants $c_i>0$, $k_i>0$.  We assume $\ell \in \mathbb{N}\backslash \left\{ 1\right\}$ as the case $\ell=1$ has been considered elsewhere.   We have
  $$\Phi'\left( z\right) = \sum_{i=1}^{\ell}  c_i k_i\, z^{k_i-1} \mbox{ and } \Phi''\left( z\right) = \sum_{i=1}^{\ell}  c_i k_i\left( k_i-1\right) \, z^{k_i-2} \mbox{.}$$
  Conditions \ref{T:Phiconditions} a) and b) are clearly satisfied.  For $\Phi''>0$ it is sufficient to have $k_i\geq 1$ for all $i$.  Since
  $$z\, \Phi''\left( z\right) = \sum_{i=1}^{\ell}  c_i k_i\left( k_i-1\right) \, z^{k_i-1} \mbox{,}$$
  Condition \ref{T:Phiconditions} d) holds if we choose $c= \max_i \left( k_i -1\right)$.\\
  
  \item $\Phi_p\left( z \right) = \ln\left( 1 + z\right) + z^p$ for any $p$ from about $1.22$ to $2$ inclusive.\\
  
  \item $\Phi_{z_0}\left( z\right) = \left( z + z_0\right) \left[ \ln\left( z + z_0\right) - 1 \right] + z_0\left( 1 - \ln z_0 \right)$ for any $z_0 \geq 1$.\\
\end{enumerate}

  \item Because $F$ is degree-one homogeneous, equation \eqref{E:sse} may be rewritten as
  \begin{equation} \label{E:sse2}
  \left< X\left( x, 0 \right), \nu\left( x, 0 \right) \right> = -\frac{1}{\psi'\left( t\right)} \Phi\left( \frac{1}{\psi\left( t\right)}F\left(  \mathcal{W} \left( x,0\right) \right) \right) \mbox{.}
  \end{equation}
However, since $\Phi$ is not homogeneous, we cannot bring the $\frac{1}{\psi\left( t\right)}$ factor outside $\Phi$ to separate out the time variable and solve for $\psi$ explicitly.  Throughout the article we will use the abbreviation $z:= \frac{1}{\psi\left( t\right)} F\left(  \mathcal{W} \left( x,0\right)  \right)$ to denote the argument of $\Phi$ and its derivatives where they occur.\\
\item In the case that $M_0$ is a sphere of radius $r$, equation \eqref{E:sse2} becomes
$$r = - \frac{1}{\psi'\left( t\right)} \Phi\left( \frac{1}{r\, \psi\left( t\right)} \right)$$
with $\psi\left( 0 \right) = 1$, which is an ordinary differential equation for $\psi$ that might or might not be solvable explicitly depending on the form of $\Phi$.  For example, if $\Phi\left( z\right) = z + z^3$ and $r=1$, the equation for $\psi$ is
$$\psi' = - \psi^{-1} - \psi^{-3} \mbox{;}$$
with initial condition $\psi\left( 0 \right)=1$ the solution can be written implicitly as
$$t= \frac{1}{2} \left( 1- \psi^2 \right) + \frac{1}{2} \ln\left( \frac{1+\psi^2}{2}\right) \mbox{.}$$
We observe that $\psi \rightarrow 0$ as $t\rightarrow \frac{1}{2} \left( 1 - \ln 2\right)$.
\end{enumerate}

Now we state our main theorem:

\begin{main} \label{T:main}
Let $M$ be a compact, convex, $n$-dimensional hypersurface, $n\geq 2$, without boundary.  Suppose $M$ satisfies equation \eqref{E:sse} where $F$ satisfies Conditions \ref{T:Fconds} and $\Phi$ satisfies Conditions \ref{T:Phiconditions} respectively.  Suppose one (or more) of the following hold
\begin{itemize}
  \item[\textnormal{i)}] $f$ is convex and $\Phi$ satisfies Conditions \ref{T:Phiconditions} d); 
  \item[\textnormal{ii)}] $f$ is concave and $\Phi$ satisfies Conditions \ref{T:Phiconditions} d); 
  \item[\textnormal{iii)}] $n=2$; 
  \item[\textnormal{iv)}] $M$ is axially symmetric. 
  \end{itemize}
If $M$ is sufficiently curvature-pinched, in the sense that there exists an $\varepsilon \in \left( 0, 1\right]$, depending on $F$, such that the principal curvatures satisfy
$$\kappa_i \geq \varepsilon \kappa_j$$
for all $i, j = 1, 2, \ldots, n$, then $M$ is necessarily a round sphere.  If none of i) to iv) above hold but $\Phi$ satisfies Condition \ref{T:Phiconditions} d) then the result also holds under a potentially-stronger pinching condition.
\end{main}

\noindent \textbf{Remark:} In case iii) above, the pinching condition on $M$ may be written explicitly as the requirement that
$$\frac{\kappa_{\max}}{\kappa_{\min}} \leq 1 + \frac{2\Phi'}{\Phi'' F}$$
holds everywhere on $M$.  In case iv) there is a similar requirement on $M$:
$$\frac{\kappa_{\mbox{axial}}}{\kappa_{\mbox{rotational}}} \leq 1 + \frac{2\Phi'}{\Phi'' F} \mbox{.}$$
These requirements on $M$ can be thought of as a replacement for the structure condition on $\Phi$ of the other cases.

While for general $\Phi$ these are conditions on $M$, it is easily checked that in the case $\Phi\left( z\right) = z^\alpha$ for $\alpha \geq 1$ they reduce to the pinching requirements as found in the earlier works \cite{Asurf} and \cite{MMW15} respectively.\\[8pt]

Let us now set up some notation, which is the same as that used elsewhere (egs \cites{Asurf, AMZ11, H90, M11}).  In particular, $g = \left\{ g_{ij}\right\}$, $A = \left\{ h_{ij}\right\}$ and $\mathcal{W} = \left\{ h^{i}_{\> j} \right\}$ denote respectively the metric, second fundamental form and Weingarten map of $M$.  The mean curvature of $M$ is
$$H= g^{ij}h_{ij} = h^{i}_{\> i}$$
and the norm of the second fundamental form is
$$\left| A \right|^{2} = g^{ij}g^{lm}h_{il}h_{jm} = h^{j}_{\> l}h_{j}^{\;\, l}$$
where $g^{ij}$ is the $\left( i, j\right)$-entry of the inverse of the matrix $\left( g_{ij}\right)$.  The norm of the trace-free component of the second fundamental form,
$$\left| A^{0}\right|^{2} = \left| A \right|^{2} - \frac{1}{n} H^{2} = \frac{1}{n} \sum_{i<j} \left( \kappa_{i} - \kappa_{j} \right)^{2} \mbox{,}$$
is identically equal to zero when $M$ is a sphere, and as in earlier work we will also set 
$$C = \kappa_{1}^{3} + \ldots + \kappa_{n}^{3} \mbox{,}$$
where the letter $C$ is chosen simply by convention.  Throughout this paper we sum over repeated indices from $1$ to $n$ unless otherwise indicated.  Raised indices indicate contraction with the metric.

We will denote by $\left( \dot{F}^{kl} \right)$ the matrix of first partial derivatives of $F$ with respect to the components of its argument:
$$\left. \frac{\partial}{\partial s} F\left( A+sB\right) \right|_{s=0} = \dot F^{kl} \left( A\right) B_{kl} \mbox{.}$$
Similarly for the second partial derivatives of $F$ we write
$$\left. \frac{\partial^{2}}{\partial s^{2}} F\left( A+sB\right) \right|_{s=0} = \ddot F^{kl, rs} \left( A\right) B_{kl} B_{rs} \mbox{.}$$
Throughout the article unless the argument is explicitly indicated we will always evaluate partial derivatives of $F$ at $\mathcal{W}$ and partial derivatives of $f$ at $\kappa\left( \mathcal{W}\right)$, where $\mathcal{W}$ is the Weingarten map of $M=M_0$.  We will further use the shortened notation $\dot f^{i} = \frac{\partial f}{\partial \kappa_{i}}$ and $\ddot f^{ij} = \frac{\partial^{2} f}{\partial \kappa_{i} \partial \kappa_{j}}$ where appropriate.

Now let us mention some geometric inequalities that will be needed in our analysis.  In view of homogeneity, it is not possible for $f$ to be strictly convex or strictly concave, as its Hessian has a zero eigenvalue in the radial direction.  However, we do have

\begin{Consign}[Lemma 7.12, \cite{A94}] \label{T:consign}
If $F$ satisfies Conditions \ref{T:Fconds} and is strictly convex in nonradial directions, then there exists $\underline{C}>0$ such that
$$\ddot F^{kl,rs} B_{kl} B_{rs} \geq \underline{C} \frac{\left| \nabla A\right|^2}{\left| A \right|} \mbox{.}$$
On the other hand, if $F$ satisfies Conditions \ref{T:Fconds} and is strictly concave in nonradial directions, then there exists $\overline{C} > 0$ such that
$$\ddot F^{kl,rs} B_{kl} B_{rs} \leq -\underline{C} \frac{\left| \nabla A\right|^2}{\left| A \right|} \mbox{.}$$
\end{Consign}

Proofs of the next result appear in \cite{M05}.

\begin{Fconvprops} \label{T:Fconvprops}
\mbox{}
\begin{itemize}
  \item[\textnormal{i)}] If $F$ satisfies Conditions \ref{T:Fconds} and $F$ is convex (concave) at $\mathcal{W}$, then at this $\mathcal{W}$,
  $$ \left| A\right|^{2} F - \dot{F}^{kl} h_{km} h^{m}_{\ \> l} H \leq \left( \geq \right) 0 \mbox{.}$$
  \item[\textnormal{ii)}]  If $F$ satisfies Conditions \ref{T:Fconds} and $F$ is convex at $\mathcal{W}$, then at this $\mathcal{W}$,
  $$FH- n \dot F^{kl} h_{k}^{\ m} h_{ml} \leq 0 \mbox{.}$$
\end{itemize}
\end{Fconvprops}

In the case where $F$ satisfies no second derivative condition other than boundedness we will also require the following geometric estimates for convex hypersurfaces.

\begin{curvest} \label{T:curvest}
Let $M$ be a closed, convex, $n$-dimensional hypersurface with pinched principal curvatures in the sense that, at every $p\in M$,
\begin{equation} \label{E:Hpinch}
  \kappa_{i} \geq \varepsilon \kappa_{j} \mbox{,}  
\end{equation}
for all $i, j = 1, \ldots, n$.  Then, at each $p\in M$ we also have:
\begin{itemize}
  \item[\textnormal{(i)}] $\left| A^{0} \right|^{2} \leq \left( \frac{n-1}{2} \right) \left( 1 - \varepsilon \right)^{2} H^{2}\mbox{,}$
   \item[\textnormal{(ii)}] $HC -\left( \left| A \right|^{2} \right)^{2} \geq \frac{\varepsilon^{2}}{n} H^{2} \left| A^{0} \right|^{2} \mbox{,}$
  \item[\textnormal{(iii)}] $\left| H \nabla_{i} h_{jk} - h_{jk} \nabla_{i} H \right|^{2} \geq \left( \frac{n-1}{2n^{2}} \right) \varepsilon^{2} H^{2} \left| \nabla A \right|^{2} \mbox{.}$
  \end{itemize}
  \end{curvest}
  
  The first of the above inequalities follows by straightforward calculation.  The second was proved in \cite{H84}, while the third appears  in \cite{CRS} attributed to Huisken.  Our constants in (ii) and (iii) above are different because of our different definition of $\varepsilon$.\\

We complete this section with a useful equation for degree zero functions $G\left( \mathcal{W}\right)$ of the Weingarten map of $M_0$.  It follows from \eqref{E:sse2} by straightforward calculations as in \cite{M11}, for example.  We set $\mathcal{L} = \dot F^{ij} \nabla_i \nabla_j$ where $\nabla_i$ denotes the covariant derivative in an orthonormal frame $\left\{ e_i\right\}$ on $M$.

\begin{eeqns} \label{T:eeqns}
If the evolving $M_t$ satisfies \eqref{E:sse} then we have for any $G\left( \mathcal{W}\right)$ smooth, degree-zero homogeneous function,
\begin{multline} \label{E:G}
  \mathcal{L} G = \left( \dot F^{ij} \ddot G^{kl, rs} - \dot G^{ij} \ddot F^{kl, rs} \right) \nabla_i h_{kl} \nabla_j h_{rs} - \frac{\psi\, \psi'}{\Phi'\left( z\right)} \left< X, e_k \right> \nabla^k G\\
  - \frac{\Phi''\left( z\right)}{\psi\, \Phi'\left( z\right)} \dot G^{ij} \nabla_i F \nabla_j F + \frac{\psi}{\Phi'\left( z\right)} \left[ z\, \Phi'\left( z\right) - \Phi\left( z\right) \right] \dot G^{ij} h_i^{\ k} h_{kj} \mbox{.}
  \end{multline}

  \end{eeqns}

\section{$F$ convex or concave} \label{S:Fcc}

If the convexity or concavity of $F$ is strict in nonradial directions, a simpler proof follows so we consider these cases first.  We use $G=\frac{H}{F}$ in equation \eqref{E:G} and find
\begin{multline} \label{E:HonF}
  \mathcal{L} \left( \frac{H}{F} \right) = - \frac{1}{F} \ddot F^{kl, rs} \nabla^i h_{kl} \nabla_i h_{rs} - \frac{2}{F} \dot F^{kl} \nabla_k F \nabla_l \left( \frac{H}{F} \right) -\frac{1}{F^2} \frac{\Phi''}{\Phi'\, \psi} \left( F\, g^{ij} - H\, \dot F^{ij} \right) \nabla_i F \nabla_j F\\
   -\frac{\psi\, \psi'}{\Phi'} \left< X, e_k \right> \nabla^k \left( \frac{H}{F} \right)+ \frac{\psi}{\Phi'} \left( z \, \Phi' - \Phi \right) \left( \left| A \right|^2 F - \dot F^{kl} h_k^{\ m} h_{ml} H \right) \frac{1}{F^2} \mbox{.}
 \end{multline}
 Let us fix $t$ in the interval of existence of the self-similar solution to \eqref{E:theflow} and consider \eqref{E:HonF} as an elliptic equation on $M$.  Since the $\ddot F$ term has a sign, we will be able to obtain a contradiction to the elliptic maximum principle if we can obtain the same sign on the $\Phi''$ and zero order terms.  It can be shown, however, that generally the $\Phi''$ term does not have a sign, so we will have to absorb it into the $\ddot F$ term and this will require the curvature pinching condition.  Specifically, using Condition \ref{T:Phiconditions} d) we  estimate
 $$-\frac{1}{\psi} F \left| \Phi''\left( \frac{1}{\psi} F\right) \right| \leq c\, \Phi'\left( \frac{1}{\psi} F\right) \mbox{,}$$
 so for $F$ convex we have
\begin{multline*}
  \ddot F^{kl, rs} \nabla^i h_{kl} \nabla_i h_{rs} + \frac{1}{F} \frac{\Phi''}{\Phi' \psi} \left( F\, g^{ij} - H\, \dot F^{ij} \right) \nabla_i F \nabla_j F\\
  \geq \overline{C} \frac{\left| \nabla A \right|^2}{\left| A \right|} - \frac{1}{F} \frac{\left| \Phi''\right|}{\Phi' \psi} \left| F\, g^{ij} - H\, \dot F^{ij} \right| \left| \dot F\right|^2 \left| \nabla A \right|^2 
  \geq \left( \frac{\overline{C}}{\left| A \right|} - \frac{1}{F^2} c \left| F\, g^{ij} - H\, \dot F^{ij} \right| \left| \dot F\right|^2 \right) \left| \nabla A \right|^2 \mbox{,}
  \end{multline*}
  while for $F$ concave we have
  \begin{multline*}
  \ddot F^{kl, rs} \nabla^i h_{kl} \nabla_i h_{rs} + \frac{1}{F} \frac{\Phi''}{\Phi' \psi} \left( F\, g^{ij} - H\, \dot F^{ij} \right) \nabla_i F \nabla_j F\\
  \leq -\overline{C} \frac{\left| \nabla A \right|^2}{\left| A \right|} + \frac{1}{F} \frac{\left| \Phi''\right|}{\Phi' \psi} \left| F\, g^{ij} - H\, \dot F^{ij} \right| \left| \dot F\right|^2 \left| \nabla A \right|^2 
  \leq \left( -\frac{\overline{C}}{\left| A \right|} + \frac{1}{F^2} c \left| F\, g^{ij} - H\, \dot F^{ij} \right| \left| \dot F\right|^2 \right) \left| \nabla A \right|^2 \mbox{.}
  \end{multline*}
  Noting that on a sphere $\left( F\, g^{ij} - H\, \dot F^{ij} \right)$ is identically equal to the zero matrix, in both cases we can ensure the $\ddot F$ and $\Phi''$ terms in \eqref{E:HonF} when combined have a sign provided $M$ is sufficiently close to a sphere in the sense that 
  $$\frac{\left| A \right|}{F^2} \left| F\, g^{ij} - H\, \dot F^{ij} \right| \left| \dot F \right|^2 \leq \frac{\overline{C}}{c} \mbox{.}$$
   
For the zero order term, using a Taylor expansion of $\Phi$ about $z$ and Condition \ref{T:Phiconditions}, (a),
   $$0= \Phi\left( 0\right) = \Phi\left( z\right) - z \Phi'\left( z\right) + \frac{1}{2} \Phi''\left( \tilde z\right) z^2 \mbox{,}$$
   for some $\tilde z$ between $0$ and $z$.  Thus
   $$z \Phi'\left( z\right) - \Phi\left( z \right) = \frac{1}{2} \Phi''\left( \tilde z\right)$$
   so if $\Phi''$ has a sign, then the sign of $z\, \Phi'\left( z\right) - \Phi\left( z \right)$ is the same.  Recalling now Lemma \ref{T:Fconvprops} we see that again if $\Phi$ is convex the zero order term in \eqref{E:HonF} has the same sign as the $\ddot F$ term.
   
  In view of Lemma \ref{T:consign}, in either case we are now in the position to apply the strict elliptic maximum principle as in \cite{M11} to see that $\frac{H}{F}$ is constant on $M$.  Using again the second derivative term and the estimate of Lemma \ref{T:consign} we obtain from \eqref{E:HonF} that $\frac{H}{F}$ identically constant on $M$ implies $\left| \nabla A \right| \equiv 0$ and thus $M$ is a sphere.\\
  
  In the cases where the convexity or concavity of $F$ is not strict in nonradial directions we proceed as in \cite{M11} with different functions $G$.  Functions of the form $G= \frac{Q}{F}$ for degree-one homogeneous and symmetric $Q\left( \mathcal{W}\right)= q\left( \kappa\right)$ satisfy
    \begin{multline} \label{E:Q}
  \mathcal{L} G = \left( \dot F^{ij} \ddot Q^{kl, rs} - \dot Q^{ij} \ddot F^{kl, rs} \right) \nabla_i h_{kl} \nabla_j h_{rs} - \frac{2}{F} \dot F^{ij} \nabla_i G \nabla_j F - \frac{\psi\, \psi'}{\Phi'\left( z\right)} \left< X, e_k \right> \nabla^k G\\
  - \frac{\Phi''\left( z\right)}{\psi\, \Phi'\left( z\right)} \dot G^{ij} \nabla_i F \nabla_j F + \frac{\psi}{\Phi'\left( z\right)} \left[ z\, \Phi'\left( z\right) - \Phi\left( z\right) \right] \dot G^{ij} h_i^{\ k} h_{kj} \mbox{.}
  \end{multline}

 In the case that $F$ is concave but not strictly concave set $Q= \left| A \right|$.  Then
  $$\frac{\partial q}{\partial \kappa_i} = \frac{\kappa_i}{\left| A \right|}> 0$$
  so the $\ddot F$ term in \eqref{E:Q} is positive.  Moreover $Q$ is strictly convex in nonradial directions, so since $\dot F$ is positive definite and attains a minimum on $M$, we have, adopting suitable coordinates
  $$\dot F^{ij} \ddot Q^{kl, rs} \nabla_i h_{kl} \nabla_j h_{rs} \geq \min_j \left( \min_M \frac{\partial f}{\partial \kappa_j} \right) \sum_i \ddot Q^{kl, rs} \nabla_i h_{kl} \nabla_i h_{rs} \geq \underline{C} \frac{\left| \nabla A\right|^2}{\left| A \right|}$$
  again using Lemma \ref{T:consign}.  Further, for the zero order term in \eqref{E:Q} we note that
$$ \dot G^{ij} h_i^{\ k} h_{kj} = \frac{1}{F^2\left| A \right|} \left( FC - \left| A\right|^2 \dot F^{ij} h_{i}^{\ k} h_{kj} \right) \geq 0
\mbox{,}$$
by Lemma \ref{T:Fconvprops}.  Thus the zero order term in \eqref{E:Q} is also nonnegative since $\Phi$ is convex.  Finally the $\Phi''$ term in \eqref{E:Q} can be absorbed by the $\ddot Q$ term in a similar way as before with sufficient curvature pinching.  Applying the elliptic maximum principle we conclude $\frac{Q}{F}$ is identically constant on $M$; \eqref{E:Q} then implies $\left| \nabla A\right| \equiv 0$ from which we conclude $M$ is a sphere.

The case of $F$ convex but not strictly convex may be proven in a similar way using the function $G=\frac{Q}{F^n}$ where $Q=K$.  Then $Q$ is strictly concave in nonradial directions.  The $\ddot Q$ term can therefore absorb the $\Phi''$ term by estimating as above, and in this case the zero order term in \eqref{E:Q} is
$$\frac{\psi}{\Phi'\left( z\right)} \left[ z\, \Phi'\left( z\right) - \Phi\left( z\right) \right] \dot G^{ij} h_i^{\ k} h_{kj} = \frac{\psi}{\Phi'\left( z\right)} \left[ z\, \Phi'\left( z\right) - \Phi\left( z\right) \right] \frac{K}{F^{n+1}} \left( F H - n \dot F^{kl} h_k^{\ m} h_{ml} \right) \leq 0$$
by Lemma \ref{T:Fconvprops}, ii) and the fact that $\Phi$ is convex.\hspace*{\fill}$\Box$

\section{The cases of $n=2$ and $M$ axially symmetric}

Both these cases can be handled with the same function $G$; we will point out the differences in each case where they arise.  First let us write Lemma \ref{T:eeqns} in a slightly different form: since
\begin{multline*}
  \frac{\Phi'\left( z\right)}{\psi} \left( \dot F^{ij} \ddot G^{kl, rs} - \dot G^{ij} \ddot F^{kl, rs} \right) \nabla_i h_{kl} \nabla_j h_{rs} - \frac{\Phi''\left( z\right)}{\psi\, \Phi'\left( z\right)} \dot G^{ij} \nabla_i F \nabla_j F\\
   = \left( \dot \Phi^{ij} \ddot G^{kl, rs} - \dot G^{ij} \ddot \Phi^{kl, rs} \right) \nabla_i h_{kl} \nabla_j h_{rs} \mbox{,}
   \end{multline*}
Lemma \ref{T:eeqns} may be rewritten as
\begin{equation} \label{E:G2}
  \overline{\mathcal{L}} G = \left( \dot \Phi^{ij} \ddot G^{kl, rs} - \dot G^{ij} \ddot \Phi^{kl, rs} \right) \nabla_i h_{kl} \nabla_j h_{rs} - \psi' \left< X, e_k \right> \nabla^k G + \left[ z\, \Phi'\left( z\right) - \Phi\left( z\right) \right] \dot G^{ij} h_i^{\ k} h_{kj} \mbox{,}
  \end{equation}
  where $\overline{\mathcal{L}} := \frac{\Phi'\left( z\right)}{\psi} \dot F^{ij} \nabla_i \nabla_j$.

For the axially symmetric case, let us denote by $\kappa_1$ the curvature in the axial direction, and $\kappa_2$ the rotational curvature.  We take $G= \frac{n \left| A^0 \right|^2}{H^2}$ so the corresponding function $g$ is (in either case)
 $$g\left( \kappa\right) = \frac{\left( n-1\right) \left( \kappa_1 - \kappa_2\right)^2}{\left[ \kappa_1 + \left( n-1\right) \kappa_2 \right]^2} \mbox{.}$$
 It is important to keep in mind that $g$ is symmetric in $\left( \kappa_1, \kappa_2\right)$ in the $n=2$ case, but not in the general axially symmetric case.  

Clearly $G\geq 0$ on $M$ and if $G\equiv 0$ then $M$ is umbilic and therefore a sphere.  So suppose for the sake of obtaining a contradiction that $G$ attains a positive maximum at some $p\in M$.  Since
$$\dot g^1 = \frac{2n\left( n-1\right) \kappa_2 \left( \kappa_1 - \kappa_2\right)}{H^3} \mbox{ and } \dot g^2 = \frac{2n\, \kappa_1\left( \kappa_2 - \kappa_1\right)}{H^3}$$
we have
$$\dot G^{ij} h_i^{\ k} h_{kj} = \dot g^1 \kappa_1^2 + \left( n-1\right) \dot g^2 \kappa_2^2 = \frac{2n\left( n-1\right) \kappa_1 \kappa_2}{H^3} \left( \kappa_1 -\kappa_2\right)^2 \geq 0 \mbox{.}$$
Thus since $\Phi$ is convex we will obtain a contradiction to the maximum principle providing the remaining terms on the right hand side of \eqref{E:G2} are positive.

 Since $F$ is degree-one homogeneous and $G$ is degree-zero homogeneous, using the condition $\nabla G = 0$ at a maximum we find using cancellation as in \cite{Asurf} that in the present setting  
\begin{align*}
& \left( \dot \Phi^{ij} \ddot G^{kl, rs} - \dot G^{ij} \ddot \Phi^{kl, rs} \right) \nabla_i h_{kl} \nabla_j h_{rs} \\
  &= -\frac{F \dot g^1}{\kappa_1 \kappa_2^2 \left( \kappa_2 - \kappa_1\right)}\left\{ \left[ \Phi'' F \kappa_1\left( \kappa_2 - \kappa_1 \right) + 2 \Phi' \kappa_1 \kappa_2 \right] \left( \nabla_1 h_{22} \right)^2 \right. \\
  & \qquad \left. + \left[ -\frac{\Phi'' F \kappa_2}{\left( n-1\right)} \left( \kappa_2 - \kappa_1\right) + 2 \Phi' \kappa_1 \kappa_2 \right]  \left( \nabla_2 h_{11} \right)^2 \right\}\\
  &=\frac{2n\left( n-1\right)F}{\kappa_1 \kappa_2 H^3}\left\{ \kappa_1 \left[ \Phi'' F \left( \kappa_2 - \kappa_1 \right) + 2 \Phi'  \kappa_2 \right] \left( \nabla_1 h_{22} \right)^2 \right. \\
  & \qquad \left. + \kappa_2 \left[ -\frac{\Phi'' F }{\left( n-1\right)} \left( \kappa_2 - \kappa_1\right) + 2 \Phi' \kappa_1 \right]  \left( \nabla_2 h_{11} \right)^2 \right\} \mbox{.}
  \end{align*}
 In view of symmetry, in the case $n=2$ we may assume $\kappa_2 \geq \kappa_1$ and for the above gradient term to be positive we require
 $$2\Phi' \kappa_1 - \Phi'' F \left( \kappa_2 - \kappa_1 \right) \geq 0 \mbox{ and } 2 \Phi' \kappa_2 + \Phi'' F \left( \kappa_2 - \kappa_1 \right) \geq 0 \mbox{.}$$
If it happens that $\Phi''=0$ at the maximum of $G$, then these conditions are obviously satisfied.  If not, the second condition is still clearly satisfied however the first is only true provided
\begin{equation} \label{E:r}
  r \leq 1 + \frac{2 \Phi'}{\Phi'' F} \mbox{.}
\end{equation}

 We conclude that if $M$ satisfies \eqref{E:r} we have a contradiction to the elliptic maximum principle unless $G$ is identically constant.  In the $n=2$ case, any surface has an umbilic point, so there is a point where $G=0$ and thus $G\equiv 0$.  It follows that $M$ is a sphere.
  
In the axially symmetric case, $\nabla_2 h_{11} \equiv 0$ (see, for example, \cite[Lemma 3.2]{MMW15}) but there is no symmetry in $\kappa_1$ and $\kappa_2$.  We require that at the maximum of $G$,  
  $$2 \Phi' \kappa_2 + \Phi'' F \left( \kappa_2 - \kappa_1 \right) \geq 0 \mbox{.}$$
  If it happens that $\Phi'' = 0$ at the maximum of $G$, then the above is clearly satisfied.  It is also clearly satisfied if $\kappa_2 \geq \kappa_1$.  Otherwise, the above becomes the requirement
$$\frac{\kappa_1}{\kappa_2} \leq 1+ \frac{2 \Phi'}{\Phi'' F} \mbox{.}$$
Thus if $M$ is axially symmetric and satisfies 
$$\frac{\kappa_{\mbox{axial}}}{\kappa_{\mbox{rotational}}} \leq 1 + \frac{2 \Phi'}{\Phi'' F}$$
 we obtain a contradiction to the elliptic maximum principle unless $G$ is identically constant.  In that case we have from \eqref{E:G2}
\begin{equation*}
  0 \equiv \frac{2n\left( n-1\right)F}{\kappa_1 \kappa_2 H^3}\left\{ \kappa_1 \left[ \Phi'' F \left( \kappa_2 - \kappa_1 \right) + 2 \Phi'  \kappa_2 \right] \left( \nabla_1 h_{22} \right)^2 \right\} + \left[ z\, \Phi'\left( z\right) - \Phi\left( z\right) \right] \dot G^{ij} h_i^{\ k} h_{kj} \mbox{.}
 \end{equation*}

Given the curvature pinching and convexity of $\Phi$, this is the sum nonpositive terms, so each must be identically equal to zero.  Since $\Phi'\, F - \Phi >0$ it must be the case that $\kappa_1 \equiv \kappa_2$ everywhere on $M$ thus all principle curvatures are equal and $M$ is a sphere.\hspace*{\fill}$\Box$
\mbox{}\\[8pt]

\section{The case of no convexity condition on $F$ but strong curvature pinching}
  In this case we take $G= \frac{H^2}{\left| A \right|^2}$ and note that $1 \leq G\leq n$.  On a sphere $G\equiv n$ so we will use the elliptic maximum principle to show it is not possible for $G$ to attain a minimum unless $G$ is constant.  We find find by a similar calculation as in \cite{M11} using Lemma \ref{T:eeqns} that
  \begin{align} \label{E:strong}
   & \frac{\left| A\right|^4}{2} \overline{\mathcal{L}} \left( \frac{H^2}{\left| A \right|^2} \right) \\ \nonumber
    &= - H \left( \left| A\right|^2 g^{ij} - H\, h^{ij} \right) \ddot \Phi^{kl, rs} \nabla_i h_{kl} \nabla_j h_{rs} + \frac{\left| A\right|^4}{H} \dot \Phi^{ij} \nabla_i H \nabla_j \left( \frac{H^2}{\left| A\right|^2} \right)\\ \nonumber
    & \quad -  \frac{\left| A\right|^4}{2} \frac{\psi' \Phi'}{\Phi} \left< X, e_k \right> \nabla^k  \left( \frac{H^2}{\left| A \right|^2} \right)
    - \dot \Phi^{ij} \left( H \nabla_i h_{kl} - h_{kl} \nabla_i H\right)  \left( H \nabla^i h_{kl} - h_{kl} \nabla^i H\right) \\
    &\quad - H  \left[ z\, \Phi'\left( z\right) - \Phi\left( z\right) \right] \left[ H\, C - \left( \left| A \right|^2\right)^2 \right] \mbox{.} \nonumber
    \end{align}
 Using Lemma \ref{T:curvest} (ii) and the fact that $\Phi$ is convex we see that the above zero order term is nonpositive.  We also have a good negative norm-like term above that we use to absorb the un-signed $\ddot \Phi$ term.  Specifically, 
\begin{align*}
 &- H \left( \left| A\right|^2 g^{ij} - H\, h^{ij} \right) \ddot \Phi^{kl, rs} \nabla_i h_{kl} \nabla_j h_{rs} \\
& = -H \left( \left| A\right|^2 g^{ij} - H\, h^{ij} \right) \frac{\Phi'\left( z\right)}{\psi} \ddot F^{kl, rs}\nabla_i h_{kl} \nabla_j h_{rs}  - H  \left( \left| A\right|^2 g^{ij} - H\, h^{ij} \right) \frac{\Phi''\left( z\right)}{\psi^2} \nabla_i F \nabla_j F\\
&\leq \sqrt{n} H \left| A^0\right| \left[ M_2\left( \varepsilon\right) + c\, M_0\left( \varepsilon\right) \right] \frac{\Phi'\left( z\right)}{\psi} \left| \nabla A\right|^2
 \end{align*}
 and in view of Lemma \ref{T:curvest} (iii) we estimate
 \begin{multline*}
 \dot \Phi^{ij} \left( H \nabla_i h_{kl} - h_{kl} \nabla_i H\right)  \left( H \nabla^i h_{kl} - h_{kl} \nabla^i H\right) \\
 \geq \frac{\Phi'\left( z\right)}{\psi} M_1\left( \varepsilon \right) \left| H \nabla_i h_{kl}- h_{kl} \nabla_i H \right| \geq \left( \frac{n-1}{2n^2} \right) \varepsilon^2 H^2 M_1\left( \varepsilon\right) \left| \nabla A \right|^2 \mbox{,}
 \end{multline*}
 where, for the curvature pinching cone
 $$K_\varepsilon = \left\{ \kappa\in \Gamma^+: \kappa_i \geq \kappa_j \mbox{ for all } 1 \leq i, j \leq n \right\}$$
 we set
 $$M_0\left( \varepsilon \right) = \sup \left\{ \frac{\left| A \right|}{F}\left| \dot F\right|^2: \kappa\in K_\varepsilon, \left| \kappa \right| = 1 \right\} \mbox{,}$$
 $$M_1\left( \varepsilon \right) = \inf\left\{ \frac{\partial f}{\partial \kappa_i} \left( \kappa \right): 1\leq i\leq n,  \kappa\in K_\varepsilon, \left| \kappa \right| = 1 \right\}$$
 and
 $$M_2\left( \varepsilon \right) =\sup \left\{ \left| D^2 f\left( \kappa\right) \left( \xi, \xi \right) \right|: \kappa\in K_\varepsilon, \left| \kappa \right| = 1, \left| \xi\right|=1 \right\} \mbox{,}$$
 all of which are positive and finite as taken over compact sets.  Observe in particular $M_1\left( \varepsilon \right)$ is attained as a positive minimum and nondecreasing in $\varepsilon$; $M_2$ and $M_0$ are nonincreasing in $\varepsilon$.
 
 Using now \ref{T:curvest} (i), from \eqref{E:strong} and the above we obtain
  \begin{equation*} 
   \frac{\left| A\right|^4}{2} \overline{\mathcal{L}} \left( \frac{H^2}{\left| A \right|^2} \right) 
    \leq Q\left( \varepsilon \right) H^2 \left|  \nabla A\right|^2+ \frac{\left| A\right|^4}{H} \dot \Phi^{ij} \nabla_i H \nabla_j \left( \frac{H^2}{\left| A\right|^2} \right)  -  \frac{\left| A\right|^4}{2} \frac{\psi' \Phi'}{\Phi} \left< X, e_k \right> \nabla^k  \left( \frac{H^2}{\left| A \right|^2} \right) \mbox{,}
    \end{equation*}
    where
    $$Q\left( \varepsilon \right) = \sqrt{\frac{n\left( n-1\right)}{2}} \left[ \left( 1-\varepsilon \right) M_2\left( \varepsilon \right) + c\, M_0\left( \varepsilon \right) \right] - \left( \frac{n-1}{2n^2} \right) \varepsilon^2 M_1\left( \varepsilon \right) \mbox{.}$$
    A direct calculation shows that $Q$ is nonincreasing in $\varepsilon$, positive for small $\varepsilon$ (weaker pinching) and negative for $\varepsilon$ closer to $1$ (strong pinching).  It follows that there is a weakest curvature pinching ratio $\varepsilon$ such that, $M$ satisfying this pinching ensures the above $\left| \nabla A \right|$ term is nonpositive and thus $\frac{H^2}{\left| A \right|^2}$ cannot attain a minimum unless it is identically constant.  In that case we get that $M$ has $\left| \nabla A\right| \equiv 0$ and thus $M$ is a sphere.\hspace*{\fill}$\Box$
 

\end{document}